# On transiso graphs of groups of order less than 32

## L. K. Mishra[1] and B. K. Sharma

Department of Mathematics, University of Allahabad
Allahabad – 211002, India.
Email: lkmp02@gmail.com, brajeshsharma72@gmail.com

**Abstract:** For a finite group $G$ and a divisor $d$ of $|G|$, the transiso graph $\Gamma_d(G)$ is a graph whose vertices are subgroups of $G$ of order $d$ and two distinct vertices $H_1$ and $H_2$ are adjacent if and only if there exist normalized right transversals $S_1$ and $S_2$ of $H_1$ and $H_2$ respectively in $G$ such that $S_1 \cong S_2$ with respect to the right loop structure induced on them. In the present paper, we have determined some finite groups $G$ for which the graphs $\Gamma_d(G)$ are complete for each divisor $d$ of $|G|$. We have also discussed the completeness of transiso graphs for groups of order less than $32$.

# 1. Introduction

Let $G$ be a finite group and $H$ be a subgroup of $G$. A normalized right transversal (NRT) $S$ of $H$ in $G$ is a subset of $G$ obtained by selecting one and only one element from each right coset of $H$ in $G$ and $1 \in S$. An NRT $S$ has an induced binary operation $\circ$ given by $\{x \circ y\} = S \cap Hxy$, with respect to which $S$ is a right loop with identity $1$ (see [9, Proposition 2.2, p.42], [7]). Conversely, every right loop can be embedded as an NRT in a group with some universal property (see [7, Theorem 3.4, p.76]). Let $\langle S \rangle$ be the subgroup of $G$ generated by $S$ and $H_S$ be the subgroup $H \cap \langle S \rangle$. Then, $H_S = \langle \{xy(x \circ y)^{-1} | x, y \in S\} \rangle$ and $H_S S = \langle S \rangle$. Identifying $S$ with the set $H \backslash G$ of all right cosets of $H$ in $G$, we get a transitive permutation representation $\chi_S : G \to Sym(S)$ defined by $\{\chi_S(g)(x)\} = S \cap Hxg$, $g \in G$, $x \in S$. The kernel $\ker \chi_S$ of this action is $Core_G(H)$, the core of $H$ in $G$. The group $G_S = \chi_S(H_S)$ is known as the group torsion of the right loop $S$ (see [7, Definition 3.1, p. 75]) which depends only on the right loop structure $\circ$ on $S$ and not on the subgroup $H$. Since $\chi_S$ is injective on $S$ and if we identify $S$ with $\chi_S(S)$, then $\chi_S(\langle S \rangle) = G_S S$ which also depends only on the right loop $S$ and $S$ is an NRT of $G_S$ in $G_S S$. One can also verify that $\ker(\chi_S|_{H_S S}: H_S S \to G_S S) = \ker(\chi_S|_{H_S}: H_S \to G_S) = Core_{H_S S}(H_S)$ and $\chi_S|_S = I_S$, the identity map on $S$. If $H$ is a corefree subgroup of $G$, then there exists an NRT $T$ of $H$ in $G$ which generates $G$ (see [2]). In this case, $G = H_T T \cong G_T T$ and $H = H_T \cong G_T$. Also $(S, \circ)$ is a group if and only if $G_S$ is trivial. Let $\mathcal{T}(G, H)$ denote the set of all normalized right transversals (NRTs) of $H$ in $G$. Two NRTs $S, T \in \mathcal{T}(G, H)$ are said to be isomorphic (denoted by $S \cong T$), if their induced right loop structures are isomorphic. A subgroup $H$ is normal in $G$ if and only if all NRTs of $H$ in $G$ are isomorphic to the quotient group $G/H$ (see [7]).

---

[1] Author is supported by University Grants Commission, India.

Throughout the paper, we will assume that $G$ is a finite group and $d$ is a divisor of the order $|G|$ of the group $G$. Let $V_d(G)$ be the set of all subgroups of $G$ of order $d$. We define a graph $\Gamma_d(G) = (V_d(G), E_d(G))$ with $\{H_1, H_2\} \in E_d(G)$ if and only if there exists $S_i \in \mathcal{T}(G, H_i)$ $(i=1,2)$ such that $S_1 \cong S_2$ with respect to the right loop structure induced on $S_i$. We will call this graph a transiso graph (see [6]). If $G$ has no subgroup of order $d$, then $\Gamma_d(G)$ is a null graph (a graph having empty vertex set and empty edge set). If $G$ has unique subgroup of order $d$, then $\Gamma_d(G)$ is an empty graph (a graph having empty edge set). We will denote transiso graph $\Gamma_d(G)$ by $\Gamma_d$ if there is no confusion about $G$. A group $G$ is called a t-group if $\Gamma_d(G)$ is a complete graph for each divisor $d$ of $|G|$.

In this paper, we have determined all t-groups of the order less than $32$. In the Section $2$, we have recalled some preliminary results related to transiso graph from the paper [6]. We have also discussed about the relation of adjacency and proved that the direct product of two t-groups of coprime order is a t-group. In the Section $3$, we have discussed about the transiso graphs of some non-abelian groups like dicyclic groups, quasidihedral groups and the groups of the order $pq, 4p, 2pq$ and $2p^2$ for distinct odd prime $p$ and $q$. We have classified all the t-groups of order less than $32$ in the Section $4$.

## 2. Preliminaries

We first recall some basic results from the paper [6] and prove some elementary results which will be used in the present paper.

**Proposition 2.1.** [6, Proposition 2.1] A subgroup of a group $G$ is always adjacent with its automorphic images in $\Gamma_d(G)$ for any divisor $d$ of $|G|$.

**Proposition 2.2.** [6, Proposition 2.2] Let $H_1$ and $H_2$ be corefree subgroups of $G$. Let $S_i \in \mathcal{T}(G, H_i)$ $(i=1,2)$ such that $S_1 \cong S_2$ and $\langle S_i \rangle = G$. Then, an isomorphism between $S_1$ and $S_2$ can be extended to an automorphism of $G$ which sends $H_1$ onto $H_2$.

**Proposition 2.3.** [6, Proposition 2.3] A finite abelian group $G$ is a t-group if and only if each Sylow subgroup of $G$ is either elementary abelian or cyclic.

**Corollary 2.4.** [6, Corollary 2.2] An elementary abelian group is a t-group.

**Proposition 2.5.** [6, Proposition 2.4] The dihedral group $D_{2n}$ of order $2n$ is a t-group.

One can easily observe that the number of vertices in the graph $\Gamma_d(D_{2n})$ is equal to the number of subgroups of $D_{2n}$ of order $d$ and is given by

$$|V_d(D_{2n})| = \begin{cases} 1 & \text{if } d \text{ is odd}. \\ \dfrac{2n}{d} & \text{if } d \text{ is even and does not divide } n. \\ \dfrac{2n}{d}+1 & \text{if } d \text{ is even and divides } n. \end{cases}$$

**Proposition 2.6.** [6, Proposition 3.1] Let $G$ be a non $p$-central finite $p$-group. Then, $\Gamma_d(G)$ is complete if and only if whenever $H$ is a non-normal subgroup of $G$ of order $p$, $G \cong H \ltimes K$ for some subgroup $K$ of $G$ with $G/L \cong K$ for any normal subgroup $L$ of $G$ of order $p$.

**Proposition 2.7.** [6] Let $p$ be an odd prime and $G$ be a non-abelian group. Then,
1. If the group $G$ is a t-group and $|G| = p^3$, then $G$ is of exponent $p$ (and hence $G \cong C_p^{\,2} \rtimes C_p$).
2. If $|G| = p^4$, then $\Gamma_p(G)$ is not a complete graph.
3. If $|G| = p^5$, then $\Gamma_p(G)$ is not complete unless $\Phi(G) = Z(G) = G' \cong C_p^2$.

Let $G$ be a finite group and $d$ be a divisor of $|G|$. Let us define a relation $\sim_d$ on the set $V_d(G)$ of all subgroups of the group $G$ of order $d$ such that two subgroups $H_1$ and $H_2$ are related by the relation $\sim_d$ if either $H_1 = H_2$ or $H_1$ and $H_2$ are adjacent in the graph $\Gamma_d(G)$. We call this relation $\sim_d$ the relation of adjacency in the graph $\Gamma_d(G)$. It is trivial that the relation $\sim_d$ is reflexive and symmetric on $V_d(G)$.

**Proposition 2.8.** If the relation $\sim_d$ defined above is a transitive relation on $V_d(G)$, then $\Gamma_d(G)$ is either a complete graph or a disjoint union of complete graphs.

**Proof.** Assume that the relation $\sim_d$ is a transitive relation on $V_d(G)$. Then, it is an equivalence relation on $V_d(G)$ and hence it gives a partition of $V_d(G)$ and each component of this partition corresponds to a complete graph.

**Lemma 2.9.** Let $H_i$ and $K_i$ ($i = 1, 2$) be subgroups of the groups $G_i$ such that there exist NRTs $S_i \in \mathcal{T}(G, H_i)$ and $T_i \in \mathcal{T}(G, K_i)$ with $S_i \cong T_i$. Then, $S_1 \times S_2 \cong T_1 \times T_2$.

**Proof.** One can easily observe that $S_1 \times S_2 \in \mathcal{T}(G_1 \times G_2, H_1 \times H_2)$, for an element $(g_1, g_2) \in G_1 \times G_2$ can be expressed as $(g_1, g_2) = (h_1 s_1, h_2 s_2) = (h_1, h_2)(s_1, s_2)$, where $h_i \in H_i$ and $s_i \in S_i$ ($i = 1, 2$). Similarly $T_1 \times T_2 \in \mathcal{T}(G_1 \times G_2, K_1 \times K_2)$. Then, the map $f \times g : S_1 \times S_2 \to T_1 \times T_2$ given by $(s_1, s_2) \in (f(s_1), g(s_2))$, is a right loop isomorphism where $f : S_1 \to T_1$ and $g : S_2 \to T_2$ are right loop isomorphisms.

**Proposition 2.10.** The direct product of two t-groups of coprime order is a t-group.

**Proof.** Let $G_1$ and $G_2$ be two t-groups of coprime order. Let $G = G_1 \times G_2$ and $H, K$ be subgroups of $G$ of same order. Then by [10, Corollary, p. 141], $H = H_1 \times H_2$ and $K = K_1 \times K_2$ for some subgroups $H_1, K_1 \in G_1$ and $H_2, K_2 \in G_2$ such that $|H_1| = |K_1| = d_1$ and $|H_2| = |K_2| = d_2$. Since $G_1$ and $G_2$ are t-groups, $H_1 \sim_{d_1} K_1$ and $H_2 \sim_{d_2} K_2$. Therefore by Lemma 2.9, the subgroups $H$ and $K$ are adjacent in the corresponding transiso graph. Hence the group $G$ is also a t-group.

**Lemma 2.11.** Let $G$ be a finite group and $H$ be a non-normal subgroup of prime order. Then, an NRT $S$ of $H$ in $G$ is either a subgroup of $G$ or $H = H_S \cong G_S$.

**Proof.** Let $S$ be an NRT of $H$ in $G$. Then, either $H_S = \{1\}$ or $H_S = H$. If $H_S = \{1\}$, then $S$ is a subgroup of $G$. Now, assume that $H_S = H$. Since $H$ is core-free, $G_S \cong H_S$. We also observe that $S$ is not a group in this case.

## 3. Transiso Graphs for some Non-abelian Groups

In this section, we have determined transiso graphs for some non-abelian groups like dicyclic groups, quasidihedral groups and the groups of the order $pq, 4p, 2pq$ and $2p^2$ for distinct odd primes $p$ and $q$.

The dicyclic group (or binary dihedral group) $Q_{4n} = \langle a, b \mid a^{2n}, a^n b^2, abab^{-1} \rangle$ is a group of order $4n$ for $n \geq 1$ (see [8, p.347]). It is a non-abelian group for $n > 1$ and it is a cyclic group for $n = 1$ (that is, $Q_4 \cong C_4$). A generalized quaternion group is a special case of the dicyclic group $Q_{4n}$ when $n = 2^k$ for some positive integer $k$.
In order to prove the Proposition 3.2, we need the following elementary lemma.

**Lemma 3.1.** A subgroup of the dicyclic group $Q_{4n}$ is either cyclic or dicyclic. Moreover if $d$ is a divisor of $4n$, then

1. there is unique subgroup (namely $\left\langle a^{\frac{2n}{d}} \right\rangle$) of $Q_{4n}$ of order $d$ if $4$ does not divide $d$.

2. there are $i$ subgroups ($\langle a^i, a^j b \rangle, 0 \leq j < i$) of order $d$ conjugate to each other if $4$ divides $d$ and $i = \dfrac{4n}{d}$ is odd.

3. a subgroup of order $d$ is either $\langle a^i \rangle$ or conjugate to one of $\langle a^i, b \rangle$ or $\langle a^i, ab \rangle$ if $4$ divides $d$ and $i = \dfrac{4n}{d}$ is even.

**Proof.** Let $H$ be a nontrivial proper subgroup of $Q_{4n}$ of order $d$. Clearly $\langle a \rangle$ is maximal cyclic subgroup of $Q_{4n}$ of index $2$. The composite homomorphism $H \to Q_{4n} \to Q_{4n}/\langle a \rangle$ is either trivial or onto with the kernel $H \cap \langle a \rangle = \langle a^i \rangle$ for unique divisor $i$ of $2n$. If the homomorphism is trivial, then $H \cap \langle a \rangle = \langle a^i \rangle$ for unique divisor $i = \dfrac{4n}{d}$ of $2n$. Therefore the subgroup $H$ is cyclic in this case.
Now, if the homomorphism is onto, then $H/\langle a^i \rangle \cong Q_{4n}/\langle a \rangle \cong C_2$. Since $H \not\subset \langle a \rangle$, $H$ has an element $a^j b$ and $a^n \subseteq \langle a^i \rangle$ for $(a^j b)^2 = a^n \in H$. Therefore $H \cap \langle a \rangle = \langle a^i \rangle$ for unique divisor

$i = \dfrac{4n}{d}$ of $n$. Now, we have an appropriate element $a^j b \in H \setminus \langle a \rangle$ where $0 \leq j < i$, such that $H = \langle a^i, a^j b \rangle$. Clearly $H$ is a dicyclic group (precisely $H \cong Q_{4 \cdot \frac{n}{i}}$) for $(a^i)^{\frac{d}{2}} = 1$, $(a^i)^{\frac{d}{4}} = (a^j b)^2$ and $(a^j b) a^i (a^j b)^{-1} = (a^i)^{-1}$.

Now, we prove the next part of the lemma.

Let H be a subgroup of $Q_{4n}$ of order $d$ and $i = \dfrac{4n}{d}$. If $d$ is not a multiple of $4$, then there is no subgroup of $Q_{4n}$ of order $d$ which is dicyclic and so $H = \langle a^{\frac{i}{2}} \rangle$ is a cyclic subgroup. If $d$ is a multiple of $4$, then there are two cases.

If $d \nmid 2n$ i.e. $i$ is odd, then $H$ can not be contained in $\langle a \rangle$ so $H$ is dicyclic subgroup of the form $\langle a^i, a^j b \rangle$. If $i \leq j$, then we can find $l$ such that $0 \leq l < i$ and $H = \langle a^i, a^l b \rangle$. Thus we conclude that $0 \leq j < i$ and hence there are $i$ subgroups of order $d$ which are conjugates.

If $d \mid 2n$ i.e. $i$ is even, then $H$ is either $\langle a^{\frac{i}{2}} \rangle$ or of the form $\langle a^i, a^j b \rangle$. Using above arguments, we see that there are $\dfrac{i}{2}$ subgroups conjugate to $\langle a^i, b \rangle$ and $\dfrac{i}{2}$ subgroups conjugate to $\langle a^i, ab \rangle$.

One can easily observe that an abelian normal subgroup of the group $Q_{4n}$ is cyclic subgroup contained in the maximal cyclic subgroup and a non-abelian normal subgroup of $Q_{4n}$ has index less than or equal to $2$.

**Proposition 3.2.** The dicyclic group $Q_{4n} = \langle a, b \mid a^{2n}, a^n b^2, abab^{-1} \rangle$ of order $4n$ is a t-group.

**Proof.** Let $d$ be a divisor of $4n$ and $i = \dfrac{4n}{d}$.

First assume that $4 \nmid d$. Then by Lemma 3.1, there is unique subgroup of $Q_{4n}$ of order $d$ and so $\Gamma_d(Q_{4n})$ is trivially a complete graph.

Now assume that $4 \mid d$ and $i$ is odd. Then by Lemma 3.1, there are $i$ subgroups of order $d$ conjugate to $\langle a^i, b \rangle$ and so $\Gamma_d(Q_{4n})$ is a complete graph.

Finally assume that $4 \mid d$ and $i$ is even. Then, a subgroup of order $d$ is either $H_1 = \langle a^{\frac{i}{2}} \rangle$ or conjugate to exactly one of $H_2 = \langle a^i, b \rangle$ or $H_3 = \langle a^i, ab \rangle$. Note that $H_1$ is a normal subgroup of $Q_{4n}$ and so its all NRTs are isomorphic to $Q_{4n}/H_1 (\cong D_{2 \cdot \frac{i}{2}})$.

Now, choose $S_2 = \left\{ a^{2j+k} b^k \mid 0 \leq j < \dfrac{i}{2}, k = 0,1 \right\}$ in $\mathcal{T}(Q_{4n}, H_2)$ and $S_3 = \left\{ a^{2j} b^k \mid 0 \leq j < \dfrac{i}{2}, k = 0,1 \right\}$ in $\mathcal{T}(Q_{4n}, H_3)$. Note that $\langle S_2 \rangle = \langle a^2, ab \rangle$ and $\langle S_3 \rangle = \langle a^2, b \rangle$. Then, $H_{S_2} = \langle S_2 \rangle \cap H_2 = \langle a^i \rangle \triangleleft \langle S_2 \rangle$ and $H_{S_3} = \langle S_3 \rangle \cap H_3 = \langle a^i \rangle \triangleleft \langle S_3 \rangle$. Therefore $G_{S_2} = G_{S_3} = \{1\}$ and hence $S_2$ and $S_3$ are groups.

Let $\circ_2$ denote the induced binary operation on $S_2$ as described in the Section 1. One can observe that, $(a^2)^{\frac{i}{2}} = (ab)^2 = (ab \circ_2 a^2)^2 = 1$. This implies that $S_2 \cong D_{2 \cdot \frac{i}{2}}$. One can similarly observe that $S_3 \cong D_{2 \cdot \frac{i}{2}}$. This shows that the graph $\Gamma_d(Q_{4n})$ is complete.

It follows from the Lemma 3.1 that the number of vertices in the graph $\Gamma_d(Q_{4n})$ is given by

$$|V_d(Q_{4n})| = \begin{cases} 1 & \text{if } 4 \text{ does not divide } d. \\ \dfrac{4n}{d} & \text{if } 4 \text{ divides } d \text{ and } \dfrac{4n}{d} \text{ is odd.} \\ \dfrac{4n}{d} + 1 & \text{if } 4 \text{ divides } d \text{ and } \dfrac{4n}{d} \text{ is even.} \end{cases}$$

The quasidihedral (or semidihedral) group $QD_{2^n} = \langle a, b \mid a^{2^{n-1}}, b^2, baba^{2^{n-2}+1} \rangle$ is a non-abelian group of order $2^n$ where $n > 4$ (see [5, p.191]). Its subgroup structure can be given by the following lemma.

**Lemma 3.3.** A proper nontrivial subgroup of the quasidihedral group $QD_{2^n}$ is either cyclic or dihedral or generalized quaternion.

**Proof.** The proof is similar to that of the Lemma 3.1.

From [5, Theorem 4.10, p.199], it follows that an abelian normal subgroup of the quasidihedral group $QD_{2^n}$ of order $d = 2^m$ is cyclic (precisely $\langle a^{2^{n-m-1}} \rangle$) and a non-abelian normal subgroup of $QD_{2^n}$ has index less than or equal to $2$.

Now, we have the following proposition from which it follows that the quasidihedral group $QD_{2^n}$ is not a t-group.

**Proposition 3.4.** Let $G$ be the quasidihedral group $QD_{2^n}$ and $d = 2^m$ be a divisor of $2^n$. Then, the graph $\Gamma_d(G)$ is complete if and only if $d \neq 2$.

**Proof.** First assume that $d \neq 2$. Then by Lemma 3.3, a subgroup of $G$ of order $d = 2^m$ is either $H_1 = \langle a^{2^{n-m-1}} \rangle \cong C_{2^m}$ or conjugate to exactly one of $H_2 = \langle a^{2^{n-m}}, b \rangle$ or $H_3 = \langle a^{2^{n-m}}, ab \rangle$. Note that $H_1$ is a normal subgroup of $QD_{2^n}$ and so its all NRTs are isomorphic to $QD_{2^n}/H_1 (\cong D_{2^{n-m}})$.

Now choose $S_2 = \{a^{2j+k}b^k \mid 0 \leq j < 2^{n-m-1}, k = 0,1\}$ in $\mathcal{T}(QD_{2^n}, H_2)$ and $S_3 = \{a^{2j}b^k \mid 0 \leq j < 2^{n-m-1}, k = 0,1\}$ in $\mathcal{T}(QD_{2^n}, H_2)$. Note that $\langle S_2 \rangle = \langle a^2, ab \rangle$ and $\langle S_3 \rangle = \langle a^2, b \rangle$. Then, $H_{S_2} = \langle S_2 \rangle \cap H_2 = \langle a^{2^{n-m}} \rangle \triangleleft \langle S_2 \rangle$ and $H_{S_3} = \langle S_3 \rangle \cap H_3 = \langle a^{2^{n-m}} \rangle \triangleleft \langle S_3 \rangle$. Therefore $G_{S_2} = G_{S_3} = \{1\}$ and hence $S_2$ and $S_3$ are groups.

Let $\circ_2$ denote the induced binary operation on $S_2$ as described in the Section 1. One can observe that, $(a^2)^{2^{n-m-1}} = (ab)^2 = (ab \circ_2 a^2)^2 = 1$. This implies that $S_2 \cong D_{2^{n-m}}$. One can similarly observe that $S_3 \cong D_{2^{n-m}}$. This shows that the graph $\Gamma_d(QD_{2^n})$ is complete.

Finally assume that $d = 2$. Then, a subgroup of $G$ of order $2$ is either $H_1 = \langle a^{2^{n-2}} \rangle$ or a conjugate to $H_2 = \langle b \rangle$. Since $H_1 \triangleleft G$, every NRT of $H_1$ in $G$ is isomorphic to $G/H_1 \cong D_{2^{n-1}}$.

Let $H$ be a non-normal subgroup of $QD_{2^n}$ of order $2$. Then, $H$ is contained in $\langle a^2, b \rangle \cong D_{2^{n-1}}$ and $H$ is a conjugate to the subgroup $\langle b \rangle$. Clearly the core $Core_G(H)$ of $H$ in $QD_{2^n}$ is trivial. Now let $S$ be an NRT of $H$ in $QD_{2^n}$. Then, the order of $H_S = H \cap \langle S \rangle$ is less than or equal to $2$. If $|H_S| = 1$, then $S = \langle S \rangle$ is a subgroup of $QD_{2^n}$. Therefore $S$ is equal to either $\langle a \rangle$ or $\langle a^2, ab \rangle \cong Q_{2^{n-1}}$.

Finally if $|H_S| = 2$, then $H_S = H$ and $\langle S \rangle = G$. Therefore $G_S \cong H_S/Core_{H_S S}(H_S) = H/Core_G(H) \cong H$. Since $G_S$ is nontrivial, $S$ is not a group. Hence $S \not\cong D_{2^{n-1}}$.

It can be trivially observed that the number of vertices in the graph $\Gamma_d(QD_{2^n})$ is equal to the number of subgroups of $QD_{2^n}$ of order $d$ and is given by

$$|V_d(QD_{2^n})| = \begin{cases} 1 & \text{if } d = 1 \text{ or } d = 2^n. \\ 2^{n-2} + 1 & \text{if } d = 2. \\ 2^{n-m} + 1 & \text{if } d = 2^m \text{ with } 0 < m < n. \end{cases}$$

**Proposition 3.5.** Let $p$ and $q$ be distinct odd primes. Then, a group of order either $pq$ or $4p$ or $2pq$ is t-group.

**Proof.** Observe that a nontrivial proper subgroup of a group of order $pq$ is a Sylow subgroup. Hence any two subgroups of same order are adjacent in corresponding transiso graph.

By classification of groups of order $4p$ (see [1, p.132-137]), a non-abelian group of order $4p$ is isomorphic to exactly one of $D_{4n}, Q_{4n}$, the alternating group $Alt(4)$ (for $p = 3$), $C_p \rtimes C_4$ (for $p \equiv 1 \mod 4$). The groups $D_{4n}$ and $Q_{4n}$ are t-groups from the propositions 2.5 and 3.2. Since any two subgroups of the group $Alt(4)$ of equal order are conjugate therefore the group $Alt(4)$ is also a t-group.

Let $H_1$ and $H_2$ be two distinct subgroups of $C_p \rtimes C_4$ of order $2$. Then, there exist unique Sylow $2$-subgroup $K_i$ of $C_p \rtimes C_4$ containing $H_i$ where $i = 1, 2$. Since $K_1$ and $K_2$ are conjugate, the subgroups $H_1$ and $H_2$ are conjugate. So $H_1$ and $H_2$ are adjacent in $\Gamma_2(C_p \rtimes C_4)$.

A non-abelian group of order $2pq$ is isomorphic to exactly one of the groups $D_{2pq}, D_{2q} \times C_p, D_{2p} \times C_q$ and $C_2 \times (C_q \rtimes C_p), (C_q \rtimes C_p) \rtimes C_2$ (when $p$ divides $q-1$) (see [4, p. 50]). $D_{2q} \times C_p, D_{2p} \times C_q$ and $C_2 \times (C_q \rtimes C_p)$ are t-groups due to the Proposition 2.10. Order of the normalizer $N_G(H)$ of a Sylow $p$-subgroup $H$ of $(C_q \rtimes C_p) \rtimes C_2$ is $2p$ and $H$ is unique Sylow $p$-subgroup of $N_G(H)$. Since all Sylow $p$-subgroups are conjugate, therefore their normalizers are also conjugate.

**Proposition 3.6.** Let $G$ be a non-abelian group of order $2p^2$ for some odd prime $p$. Then, the group $G$ is t-group if and only if $G$ is isomorphic to either the dihedral group $D_{2p^2}$ or $(C_p)^2 \rtimes C_2$.

**Proof.** It is well known that a non-abelian group of order $2p^2$ is isomorphic to exactly one of the groups $D_{2p^2}$, $(C_p)^2 \rtimes C_2$ and $C_p \times D_{2p}$ (see [1, p.132-137]).

Let $G = \langle a,b,c \mid a^p, b^p, c^2, [a,b], (ac)^2, (bc)^2 \rangle \cong (C_p)^2 \rtimes C_2$. Then, all subgroups of $\langle a,b \rangle \cong (C_p)^2$ are normal in $G$ and their quotients are dihedral groups $D_{2p}$. Hence $\Gamma_p(G)$ is a complete graph. Now $\Gamma_{2p}(G)$ is also complete as there are several NRTs of a subgroup $H$ of $G$ order $2p$ which are isomorphic to the cyclic group of order $p$. So $G$ is a t-group.

Now, let $G \cong C_p \times D_{2p} = \langle a,b,c \mid a^p, b^p, c^2, [a,b], [a,c], (bc)^2 \rangle$. Then, it is obvious that $\langle a \rangle$ and $\langle b \rangle$ are normal subgroups of $G$ of order $p$ such that $G/\langle a \rangle \cong D_{2p}$ and $G/\langle b \rangle \cong C_{2p}$. Hence $\Gamma_p(G)$ is not a complete graph.

## 4. Classification of t-groups of Order less than 32

Abelian t-groups are already determined by Proposition 2.3 which tells that a finite abelian group $G$ is a t-group if and only if it is isomorphic to the direct sum of a cyclic group $C$ and a direct sum $A$ of some elementary abelian groups, where $|A|$ and $|C|$ are coprime.

Non-abelian groups of the order 12, 20, 21, 28 and 30 are t-groups by Proposition 3.5 and a non-abelian t-group of the order 18 can be determined by Proposition 3.6. By Propositions 3.2 and 2.5, it is clear that the non-abelian groups of order 8 and $2p$ (for odd prime $p \leq 13$) are t-groups. In Propositions 4.1 and 4.3, we have determined non-abelian t-groups of the order 16 and 24 respectively. We recall that a finite $p$-group $P$ is $p$-central if each subgroup of $P$ of order $p$ is contained in the center $Z(P)$.

**Proposition 4.1.** Let $G$ be a non-abelian group of order 16. Then, the group $G$ is a t-group if and only if $G$ is isomorphic to either dihedral group $D_{16}$ or dicyclic group $Q_{16}$.

**Proof.** If $G$ is a 2-central group, then it is isomorphic to one of the groups $Q_8$, $C_4 \rtimes C_4$ and $C_2 \times Q_8$ (see [11]). By Proposition 3.2, $Q_{16}$ is a t-group. The group $C_4 \rtimes C_4 = \langle a,b \mid a^4, b^4, abab^{-1} \rangle$ has three normal subgroups $\langle a^2 \rangle$, $\langle b^2 \rangle$ and $\langle a^2 b^2 \rangle$ of order 2 with quotient groups isomorphic to the groups $C_4 \times C_2$, $D_8$ and $Q_8$ respectively. Therefore the graph $\Gamma_2(C_4 \rtimes C_4)$ is not complete and hence $C_4 \rtimes C_4$ is not a t-group. The group $C_2 \times Q_8 = \langle a,b,c \mid a^2, b^4, b^2c^2, [a,b], [a,c], bcbc^{-1} \rangle$ is not a t-group, for it has three normal subgroups $\langle a \rangle$, $\langle b^2 \rangle$ and $\langle ab^2 \rangle$ of order 2 with quotient groups isomorphic to the groups $Q_8$, $(C_2)^3$ and $Q_8$ respectively. Therefore $C_4 \rtimes C_4$ is not a t-group. If $G$ is a non 2-central group which is also a t-group, then $\Gamma_2(G)$ is a complete graph and hence by Proposition 2.6, $G$ should be isomorphic to a nontrivial semidirect product $H \ltimes K$ of a non-normal subgroup $H$ of $G$ of order 2 and a normal subgroup $K$ of $G$ of order 8 such that for any normal subgroup $L$ of $G$ of order 2, $K$ is isomorphic to $G/L$. By [11], we observe that there are five groups $(C_4 \times C_2) \rtimes_1 C_2$, $C_8 \rtimes C_2$, $QD_{16} = C_8 \rtimes_1 C_2$, $D_{16} = D_8 \rtimes C_2$ and $(C_4 \times C_2) \rtimes_2 C_2$ of required semidirect product type. Proposition 2.5 asserts that the group $D_{16}$ is a t-group and the group $QD_{16}$ is not a t-group by Proposition 3.4. The groups $(C_4 \times C_2) \rtimes_1 C_2$, $C_8 \rtimes C_2$ and $(C_4 \times C_2) \rtimes_2 C_2$ have normal subgroups of order $2$ such that corresponding quotient groups are isomorphic to $D_8$, $C_4 \times C_2$ and $(C_2)^3$ respectively (see [11]). Therefore these groups are not t-groups.

**Lemma 4.2.** Let $G$ be the group $C_2 \times Alt(4)$. Then, the graph $\Gamma_2(G)$ is not a complete graph.

**Proof.** First note that $N = C_2 \times \{1\}$ is a normal subgroup of $G = C_2 \times Alt(4)$ of order $2$, where $I$ is the identity element of $Alt(4)$ and every NRT of $N$ in $G$ is isomorphic to $G/N \cong Alt(4)$.
Now, choose a non-normal subgroup $H$ of $G$ of order $2$ which is contained the subgroup $C_2 \times Alt(4)$ of $G$.
Let $S$ be an NRT of $H$ in $G$. Note that $S' = S \cap (C_2 \times Alt(4))$ is an NRT of $H$ in $C_2 \times Alt(4)$ and $\langle S' \rangle = C_2 \times Alt(4)$. Hence by Lemma 2.11, $S$ can not be a group. Thus, the subgroups $H$ and $N$ are not adjacent in the graph $\Gamma_2(G)$, that is, the graph $\Gamma_2(G)$ is not complete.

**Proposition 4.3.** Let $G$ be a non-abelian group of order $24$. Then, the group $G$ is a t-group if and only if $G$ is isomorphic to a semidirect product of two t-groups of coprime order except the groups $C_2 \times Alt(4)$ and $(C_2 \times C_6) \rtimes C_2$.

**Proof.** We know that there are $12$ non isomorphic non-abelian groups of order $24$ (see [1, p.101-104]) and $9$ of them are semidirect product of two t-groups of coprime order.
It is obvious that the groups $C_3 \rtimes C_8$ and $SL(2,3)$ are t-groups, for any two subgroups of respective groups of equal order are conjugate. The groups $Q_{24}$ and $D_{24}$ are also t-groups by Propositions 3.2 and 2.5 respectively. By Proposition 2.10, we see that the groups $C_3 \times D_8$, $C_3 \times Q_8$ and $C_2 \times D_{12} \cong (C_2)^2 \times D_6$ are t-groups. It is clear from [6, Example 2.2] that the symmetric group Sym(4) is not a t-group. One can observe that $\langle a^2 \rangle \times Alt(3) \cong C_6$ and $\{1\} \times Sym(3)$ are normal subgroups of the group $\langle a \rangle \times Sym(3) \cong C_4 \times D_6$ such that their quotient groups are $(C_2)^2$ and $C_4$ respectively. So $\Gamma_6(C_4 \times D_6)$ is not a complete graph and hence the group $C_4 \times D_6$ is not a t-group. Similarly $C_2 \times D_{12}$ is not a t-group since there are two normal subgroups $C_2 \times \{1\}$ and $\{1\} \times Z(D_{12})$ of order $2$ such that their quotient groups are $D_{12}$ and $Q_{12}$.
Now, consider $G = (C_2 \times C_6) \rtimes C_2$. It has a normal subgroup $H$ of order $2$ such that $G/H \cong D_{12}$. Let $K$ be a subgroup of $G$ of order $2$ contained in the subgroup isomorphic to $D_{12}$. Then, there is no NRT $S \in \mathcal{T}(G,H)$ such that $S = D_{12}$, for otherwise $S = \langle S \rangle$ and $S \cap H = H$ which contradicts the fact that $S$ is an NRT. Therefore the group $(C_2 \times C_6) \rtimes C_2$ is not a t-group. Finally by Lemma 2.11, the group $C_2 \times Alt(4)$ is not a t-group.

**Acknowledgment:** Authors are thankful to Prof. R. P. Shukla, Department of Mathematics, University of Allahabad, India and Dr. Vipul Kakkar, School of Mathematics, Harish-Chandra Research Institute, Allahabad, India for suggesting this problem and their valuable discussions.